# Combinatorial Problems about Free Groups and Algebras

Alexander A. Mikhalev[*]

Vladimir Shpilrain  and  Jie-Tai Yu[†]

**Abstract**

This is a survey of recent progress in several areas of combinatorial algebra. We consider combinatorial problems about free groups, polynomial algebras, free associative and Lie algebras.

Our main idea is to study automorphisms and, more generally, homomorphisms of various algebraic systems by means of their action on "very small" sets of elements, as opposed to a traditional approach of studying their action on subsystems (like subgroups, normal subgroups; subalgebras, ideals, etc.) We will show that there is a lot that can be said about a homomorphism, given its action on just a single element, if this element is "good enough". Then, we consider somewhat bigger sets of elements, like, for example, automorphic orbits, and study a variety of interesting problems arising in that framework.

One more point that we make here is that one can use similar combinatorial ideas in seemingly distant areas of algebra, like, for example, group theory and commutative algebra. In particular, we use the same language of "elementary transformations" in different contexts and show that this approach appears to be quite fruitful for all the areas involved.

[*]Partially supported by RFBR and by INTAS.
[†]Partially supported by RGC Fundable Grant 334/024/0002 and CRCG Grant 335/024/0008.

# Contents

# 1 Introduction

The purpose of this survey is to collect in one place a number of results that appeared during the last five years and gave new life to rather old areas of algebra, such as combinatorial group theory (which is over a hundred years old), combinatorial commutative algebra (which emerged in the sixties, after Gröbner bases were introduced), and combinatorial non-commutative algebra (which can be traced back to the forties).

Of course, we were not attempting to collect *all* recent results that influenced those areas; that would take several books rather than a brief survey like this one. Our focus here is on the dynamics of homomorphisms of free groups and algebras, and, in particular, on various methods of distinguishing automorphisms among arbitrary endomorphisms. Roughly speaking, we interpret "combinatorial methods" here as "methods of elementary transformations". This allows us to use basically the same language (which is essentially matrix language) and the same general ideas for seemingly distant algebraic systems, like, for instance, free groups and polynomial algebras. This shifting between different areas appears to be quite fruitful for all the areas involved. A good illustration to that might be the following: in [78] and [79], V. Shpilrain introduced, for the purpose of recognizing automorphisms, a concept of *rank* of an element of a free group or a free Lie algebra. Then, this concept got a new meaning through the work of A. A. Mikhalev and A. Zolotykh [61] (on free Lie algebras) and U. Umirbaev [93] (on free groups). This work revealed a remarkable duality between the rank of an element of a free Lie algebra or a free group, and the rank (in the "usual" sense, i.e., the minimal number of generators) of the (one-sided) ideal (of the enveloping free associative algebra or of a free group ring, respectively) generated by (non-commutative) partial derivatives of this element. Later on, V. Shpilrain and J.-T. Yu [83] adopted this duality for two-variable polynomial algebras to produce several new results and algorithms in that classical area, which have prospective applications to a classification of non-singular plane algebraic curves.

## Notation

We sometimes use different notation in different sections of this survey to better serve purposes of a particular section. To avoid confusion, we gather below some basic types of notation. Additional notation and definitions can be found in the Background section.

All groups and algebras in this survey are finitely generated. By $F_n$ we denote the free group of rank $n$; by $A_n$, or $A(X)$, or $K\langle x_1, \ldots, x_n \rangle$ the free associative algebra (of rank $n$) over a field $K$; by $L_n$ or $L(X)$ the free Lie algebra (of rank $n$), and by $P_n$ or $K[x_1, \ldots, x_n]$ the polynomial algebra in $n$ variables. Sometimes, however, when we do not want to be very specific, we write $F(X)$ for a free algebra (freely generated by a set $X$) in a variety of (associative or non-associative) algebras.

## 2 Test elements and automorphic orbits

An element $u$ is called a *test element* if for any endomorphism $\varphi$, $\varphi(u) = u$ implies $\varphi$ is an automorphism. This definition was explicitly given by V. Shpilrain in [78], although the idea goes back to J. Nielsen [69]. In this section, we have gathered new results on test elements; for an overview of old ones, we refer to [78] and [81].

### 2.1 Test elements and automorphic orbits in free groups

An element $u \in F_n$ is called a *test element* (for monomorphisms) if for any endomorphism (any monomorphism) $\varphi$, $\varphi(u) = u$ implies $\varphi$ is an automorphism.

A classical result of Nielsen [69] is: an endomorphism $x \to f$; $y \to g$ of the free group $F_2 = \text{gp}\langle x, y \rangle$ is an automorphism if and only if $[f, g]$ is conjugate to $[x, y]$ or $[x, y]^{-1}$. Hence the commutator $[x, y] = xyx^{-1}y^{-1}$ is a test element of $F_2$.

Other classes of test elements were found by Zieschang, Rosenberger, Rips, Dold and others. A more detailed information on all those results can be found in [78] and [81]. Here we just give a couple of examples of test elements in the free group $F_n = \text{gp}\langle x_1, \ldots, x_n \rangle$: $x_1^k x_2^k \ldots x_n^k$, $k \geq 2$; $[x_1, x_2] \cdot [x_3, x_4] \cdot \ldots \cdot [x_{n-1}, x_n]$ (if $n$ is even); $[x_1, \ldots, x_n]$.

Some other examples of test elements are given in an important paper by Turner [88], where he also proves the following theorem. Recall that a *retract* $R$ of a group $G$ is a subgroup with the property that there is a normal subgroup $N$ of $G$ such that $R = G/N$. Otherwise stated, a *retraction* $\rho : G \to G$ is a homomorphism such that $\rho^2 = \rho$ and a retract is the image of a retraction. Free factors are retracts but not all retracts of free groups are free factors.

**Theorem 1 (Turner [88])** *The test elements in $F_n$ are the elements not contained in proper retracts. The test elements for monomorphisms in $F_n$ are the elements not contained in proper free factors.*

We note that there is an algorithm for deciding whether or not a given element of a free group belongs to a proper retract. This algorithm, which is due to Comerford [12], is based on results of Makanin and Razborov on solving equations in free groups.

The paper of Comerford [12] contains additional examples of test elements. Yet further examples can be found in a recent paper by Fine, Isermann, Rosenberger, and Spellman [26].

We also mention some results on test elements in other groups. Durnev [23] established the Nielsen commutator test for the free metabelian group of rank 2. On the other hand, it may come as a surprise that the Nielsen commutator test is not valid for free solvable non-metabelian groups of rank 2 – see [29].

The idea used in the definition of a test element leads to another interesting and important problem: how many elements are needed to completely determine an endomorphism of the free group $F_n$? More formally, we look for elements $g_1, \ldots, g_k \in F_n$ with the following property: whenever $\varphi(g_i) = \psi(g_i), i = 1, \ldots, k$, for some endomorphisms $\varphi$, $\psi$, of $F_n$, it follows that $\varphi = \psi$? It is obvious that $n$ elements are enough; on the other hand, it is easy to see that one element is not enough even to determine an *automorphism* since any element $u$ is fixed by the conjugation by any power of $u$. Furthermore, if the image of an endomorphism $\varphi$ is cyclic, then $\varphi$ cannot be determined by its values on less than $n$ elements! This situation, however, seems to be exceptional; our conjecture is that endomorphisms with non-cyclic images are completely determined by their values on just two elements.

S. Ivanov [31] proved this conjecture for monomorphisms:

**Theorem 2 (S. Ivanov [31])** *There are elements $u, v \in F_n$ with the following property: whenever $\varphi(u) = \psi(u)$ and $\varphi(v) = \psi(v)$ for some monomorphisms $\varphi$, $\psi$ of $F_n$, it follows that $\varphi = \psi$.*

We also note that if we consider natural extensions of endomorphisms of the free group $F_n$ to endomorphisms of the free group ring $\mathbf{Z}(F_n)$, then it becomes possible to completely determine any endomorphism by its value on just a single element of $\mathbf{Z}(F_n)$ – see [81].

Another ramification of the test element idea leads to the following problem due to Shpilrain [78]:

– Denote an orbit $\{\psi(u) \mid \psi \in Aut(F_n)\}$ by $\mathrm{Orb}(u)$. Suppose that $\varphi$ is an endomorphism of $F_n$ such that

$$\varphi(\mathrm{Orb}(h)) \subseteq \mathrm{Orb}(h)$$

for some non-identity element $h$ of $F_n$. Is then $\varphi$ an automorphism of $F_n$?

For $n = 2$, a positive solution was given independently by S. Ivanov in [32] and by V. Shpilrain in [82].

Probably the most interesting special case of this problem is where the orbit consists of primitive elements of the group $F_n$. S. Ivanov [32] answered this special case of the problem in the affirmative under an additional assumption on $\varphi$ having a *primitive pair* of elements in the image.

## 2.2 Test elements, retracts, and automorphic orbits in free algebras

In this section, we treat non-associative free algebras; for free associative algebras, see the next section.

Results of this section are due to A. A. Mikhalev and J.-T. Yu [54, 55, 56]. The proofs essentially used Theorem 35 (see Section 6.3).

In this section, $F(X)$ denotes a free $K$-algebra (without a unit element) on a finite set $X$ of free generators in one of the following varieties of algebras over a field $K$: the variety of all algebras; the variety of Lie algebras; the variety of Lie $p$-algebras; varieties of color Lie superalgebras; varieties of color Lie $p$-superalgebras; varieties of non-associative commutative and anti-commutative algebras.

We will also need the definition of *rank* of an element $u$ of a free algebra $F(X)$ (see [79]) as the minimal number of generators from $X$ on which an automorphic image of $u$ can depend.

**Theorem 3** *Let $K$ be a field, $X = \{x_1, \ldots, x_n\}$, $u \in F = F(X)$, $\operatorname{rank}(u) = n$, $\varphi$ a monomorphism of $F$ such that $\varphi(u) = u$. Then $\varphi$ is an automorphism of $F$.*

For free Lie algebras and free color Lie superalgebras Theorem 3 was proved by A. A. Mikhalev and A. A. Zolotykh [62].

**Theorem 4** *Let $X = \{x_1, \ldots, x_n\}$, $u \in F(X)$, $\operatorname{rank}(u) = n$, and let $\varphi$ be a monomorphism of the free algebra $F = F(X)$. Then $\varphi$ is an automorphism of $F$ if and only if the element $\varphi(u)$ belongs to the orbit of $u$ under the action of the automorphism group of $F$ ($\varphi(u) \in \operatorname{Orb}(u)$).*

**Theorem 5** *Let $X = \{x_1, \ldots, x_n\}$. Test elements of the free algebra $F = F(X)$ are precisely those elements not contained in any proper retract of $F$.*

**Theorem 6** *Let $K$ be a field, $X = \{x_1, \ldots, x_n\}$, and let $H$ be a nonzero subalgebra of the free algebra $F = F(X)$.*

Then $H$ is a proper retract of $F$ if and only if there exist a set $Y = \{y_1, \ldots, y_n\}$ of free generators of the algebra $F$, an integer $r$, $1 \leq r < n$, and a set $U = \{u_1, \ldots, u_r\}$ of free generators of the algebra $H$ such that

$$u_i = y_i + u_i^*, \ 1 \leq i \leq r,$$

where the elements $u_i^*$ belongs to the ideal of the algebra $F$ generated by the free generators $y_{r+1}, \ldots, y_n$.

Theorem 6 shows that retracts of free algebras $F(X)$ of finite rank have the same description as retracts of free groups of finite rank.

**Theorem 7** *Let $K$ be a field, $X = \{x_1, \ldots, x_n\}$, $u$ be a nonzero element of the algebra $F = F(X)$, $\mathrm{Orb}(u)$ the automorphic orbit of the element $u$, $\varphi$ an endomorphism of $F$ such that*

$$\varphi(\mathrm{Orb}(u)) \subseteq \mathrm{Orb}(u).$$

*Then $\varphi$ is a monomorphism of $F$.*

Let $K$ be a field, $\mathrm{char}\, K \neq 2$, $X = \{x_1, \ldots, x_n\}$ be a $G$-graded set. By $L_X$ we denote the free color Lie superalgebra $L(X)$ in the case when $\mathrm{char}\, K = 0$, and the free color Lie $p$-superalgebra $L^p(X)$ in the case when $\mathrm{char}\, K = p > 2$.

**Theorem 8** *Let $K$ be a field, $\mathrm{char}\, K \neq 2$, $X = \{x_1, \ldots, x_n\}$, and let $u$ be a nonzero element of the algebra $L = L_X$, $\mathrm{Orb}(u)$ the automorphic orbit of the element $u$, $\varphi$ an endomorphism of $L$ such that $\varphi(\mathrm{Orb}(u)) \subseteq \mathrm{Orb}(u)$. Then $\varphi$ is an automorphism of $L$.*

Note that Theorem 8 gives a positive solution of Shpilrain's problem (see the previous section) for free Lie algebras and superalgebras. In the case where $u$ is a primitive element of $L_X$, Theorem 8 was proved by A. A. Mikhalev and A. A. Zolotykh [58, 59].

## 2.3 Test elements and retracts in polynomial algebras

W. Dicks [16, 17] established the following "commutator test" for a free associative algebra $K\langle x, y\rangle$ of rank two: an endomorphism $x \to u$; $y \to v$ of the algebra $K\langle x, y\rangle$ is an automorphism if and only if $[u, v] = uv - vu = c \cdot [x, y]$ for some $c \in K^*$. Thus, in our language, the commutator $[x, y] = xy - yx$ is a test element of $K\langle x, y\rangle$.

The first example $x_1^2 + \ldots + x_n^2$ of a test polynomial in a polynomial algebra $\mathbf{R}[x_1, \ldots, x_n]$ was given by A. van den Essen and V. Shpilrain in [25]. However, in [21], V. Drensky and J.-T. Yu showed that this is not a test polynomials for $\mathbf{C}[x_1, \ldots, x_n]$. Therefore, whether or not a polynomial from $K[X]$ is a test polynomial, depends on the properties of the ground field $K$.

In that paper, several test polynomials were obtained for both $K[X]$ and $K\langle X\rangle$, a free associative algebra of finite rank. For example, the element

$$[x_1, x_2] \cdot \ldots \cdot [x_{2n-1}, x_{2n}]$$

is a test polynomial for a free associative algebra $K\langle x_1, \ldots, x_{2n}\rangle$.

Now we shall discuss *retracts* of polynomial algebras. Recall that a subalgebra $B$ of an algebra $A$ is a retract if there is an ideal $I$ of $A$ such that $A = B \bigoplus I$. The presence of other, equivalent, definitions of retracts (see e.g. [14]) provides several different methods of studying and applying them, and brings together ideas from different areas of algebra.

Let $K[x, y]$ be the polynomial algebra in two variables over a field $K$ of characteristic zero. Costa [14] proved that every proper retract of $K[x, y]$ is isomorphic to a polynomial $K$-algebra in one variable. V. Shpilrain and J.-T. Yu [85] obtained the following characterization of retracts of $K[x, y]$:

**Theorem 9 ([85])** *Let $K[p]$ be a retract of $K[x, y]$. Then there is an automorphism $\psi$ of $K[x, y]$ that takes the polynomial $p$ to $x + y \cdot q$ for some polynomial $q = q(x, y)$. A retraction for $K[\psi(p)]$ is given then by $x \to x + y \cdot q$; $y \to 0$.*

Geometrically, Theorem 9 says that (in case $K = \mathbf{C}$) every polynomial retraction of a plane is a "parallel" projection (sliding) on a fiber of a coordinate polynomial (which is isomorphic to a line) along the fibers of another polynomial (which generates a retract of $K[x, y]$).

The proof [85] of this result is based on a famous Abhyankar-Moh theorem [1].

Theorem 9 yields another useful characterization of retracts of $K[x, y]$:

**Corollary** *A polynomial $p \in K[x, y]$ generates a retract of $K[x, y]$ if and only if there is a polynomial mapping of $K[x, y]$ that takes $p$ to $x$. The "if" part is actually valid for a polynomial algebra in arbitrarily many variables.*

The following observation on retracts of a polynomial algebra in arbitrarily many variables is based on a result of Connell and Zweibel [13]:

**Proposition** *Let $R$ be a proper retract of $K[x_1, ..., x_n]$ generated by polynomials $p_1, ..., p_n$, $n \geq 2$. Then $p_1, ..., p_n$ are algebraically dependent.*

We note that it is an open problem whether or not any retract of $K[x_1, ..., x_n]$, $n \geq 3$, can be generated by algebraically independent polynomials. This problem is related to (some forms of) the well-known *cancellation problem* – see [14] for discussion.

Theorem 9 also yields a characterization of retracts of a free associative algebra $K\langle x, y \rangle$ if one uses a natural lifting:

**Theorem 10 ([85])** *Let $R$ be a proper retract of $K\langle x, y \rangle$. There is an automorphism $\psi$ of $K\langle x, y \rangle$ that takes $R$ to $K\langle v \rangle = K[v]$ for some element $v$ of the form $x + w(x, y)$, where $w(x, y)$ belongs to the ideal of $K\langle x, y \rangle$ generated by $y$.*

# 3 Primitive elements and the inverse function theorem

A notorious problem of commutative algebra (actually, it cuts across several different areas of mathematics), the Jacobian conjecture, has inspired research on similar problems in a non-commutative situation. In this section, we review some related results.

## 3.1 The Jacobian conjecture

Let $K$ be a field, char $K = 0$, $X = \{x_1, \ldots, x_n\}$, $K[X]$ the polynomial algebra in $n$ variables. For any endomorphism $\varphi$ of $K[X]$, one can consider the Jacobian matrix

$$J(\varphi) = \left( \frac{\partial \varphi(x_i)}{\partial x_j} \right), \qquad 1 \leq i, j \leq n.$$

Then we have:

**Jacobian Conjecture.** (O. H. Keller [36]) If for a polynomial mapping $\varphi$ of $K[X]$, the corresponding Jacobian matrix $J(\varphi)$ is invertible, then $\varphi$ is an automorphism of $K[X]$.

The Jacobian conjecture is still open for $n \geq 2$ (for $n = 1$ it is obviously true). For a history and background on this problem, we refer to [5]. A more recent survey is [24].

There are numerous reductions and re-formulations of the Jacobian conjecture, as well as partial results (mostly in two-variable case). Of those that are in line with the present survey, we mention a nice result of Formanek [27] that implies the following: let $\varphi$ be a polynomial mapping of $K[x, y]$ with invertible Jacobian matrix. Suppose $\varphi(K[x, y])$ contains a coordinate polynomial. Then $\varphi$ is an automorphism.

We call a polynomial $p \in K[x_1, \ldots, x_n]$ *coordinate* if it can be included in a generating set of cardinality $n$ of the algebra $K[x_1, \ldots, x_n]$.

It is interesting to note that the Jacobian conjecture for the ground field $\mathbf{C}$ is equivalent to that for the ground field $\mathbf{R}$ – see [5]. While methods of algebraic geometry are only applicable to $\mathbf{C}[X]$, combinatorial methods might give better results for $\mathbf{R}[X]$. For example, J.-T. Yu [97] reduced the Jacobian conjecture for $\mathbf{R}[X]$ to the so-called *positive case* and solved the *negative case*:

**Theorem 11** *To prove the Jacobian conjecture, one only needs to consider the case where $\varphi(x_i) = x_i + H_i^{(2)} + H_i^{(3)} + H_i^{(4)} \in \mathbf{R}[x_1, \ldots, x_n]$, where $H_i^{(j)}$ are homogeneous of degree $j$ and all coefficients in $\varphi(x_i)$ are non-negative.*

**Theorem 12** *Let $\varphi(x_i) = x_i - H_i \in \mathbf{R}[x_1, \ldots, x_n]$, where*

$$J(f_1, \ldots, f_n) \in GL_n(\mathbf{R}[x_1, \ldots, x_n]),$$

$deg(H_i) \geq 2$, *and all coefficients of $H_i$ are non-negative. Then $\varphi$ is an automorphism.*

V. Shpilrain and J.-T. Yu [85] reduced the two-variable Jacobian conjecture to the following

**Conjecture "R"**. If for a pair of polynomials $p, q \in K[x, y]$, the corresponding Jacobian matrix is invertible, then $K[p]$ is a retract of $K[x, y]$.

This statement is formally much weaker than the Jacobian conjecture since, instead of asking for $p$ to be a coordinate polynomial, we only ask for $p$ to generate a retract, and this property is much less restrictive as can be seen from Theorem 9. However, the point is that these conjectures are actually equivalent:

**Theorem 13 ([85])** *Conjecture "R" implies the Jacobian conjecture.*

We conclude this section by mentioning an analog of Shpilrain's "primitive-to-primitive" problem (see Section 2) for polynomial algebras. A. van den Essen and V. Shpilrain [25] showed that if a polynomial mapping $\varphi$ of $K[x, y]$ takes every coordinate polynomial to a coordinate one, then $\varphi$ is an automorphism of $K[x, y]$. They also showed that if the $n$-variable Jacobian conjecture is true, then the "primitive-to-primitive" problem has an affirmative answer for a polynomial algebra in $(n+1)$ variables. There was a hope for constructing a counterexample to the Jacobian conjecture by constructing a "coordinate-preserving" polynomial mapping which is not an automorphism. However, Jelonek [34] recently proved that for polynomial algebras over the field of complex numbers, every "coordinate-preserving" polynomial mapping is an automorphism.

## 3.2 Non-commutative Jacobian conjecture

**Theorem 14** *If $X = \{x_1, \ldots, x_n\}$ and $\varphi$ is an endomorphism of a free associative algebra $A(X)$, then $\varphi$ is an automorphism if and only if the Jacobian matrix*

$$\left(\frac{\partial \varphi(x_j)}{\partial x_i}\right), \qquad 1 \leq i, j \leq n,$$

*is invertible over $A(X)^e$.*

Here $A(X)^e$ is the tensor product $A(X) \otimes_K A(X)$ with the multiplication given by $(a \otimes b)(c \otimes d) = ac \otimes db$. Partial derivatives here are the components of the universal derivation – see below.

This result for $n = 2$ is due to W. Dicks and J. Lewin [19], and for arbitrary $n$ to A. H. Schofield [74].

A. A. Mikhalev and A. A. Zolotykh [63] proved a similar result for free associative algebras over an associative commutative ring.

For any element $a$ of the free associative algebra $A(X)$, we have the unique presentation in the form $a = \alpha \cdot 1 + x_1 a_1 + x_2 a_2 + \cdots$, where only finite number of elements $a_i \in A(X)$ are nonzero, and $\alpha \in K$. We call the element $a_i$ the *right Fox partial derivative of the element $a$ with respect to $x_i$*, and we use the notation $a_i = \dfrac{\partial a}{\partial x_i} = \dfrac{\partial}{\partial x_i} a$.

Thus we have the operators $\dfrac{\partial}{\partial x_i}$. These operators are linear mappings

$$\frac{\partial}{\partial x_i} \colon A(X) \to A(X)$$

such that

$$\frac{\partial}{\partial x_i}(x_j) = \delta_{ij}$$

and

$$\frac{\partial}{\partial x_i}(uv) = \frac{\partial}{\partial x_i}(u)\,v + \sigma(u)\,\frac{\partial}{\partial x_i}(v),$$

where $\sigma \colon A(X) \to A(X)$ is the homomorphism defined by $\sigma(x_i) = 0$ for all $x_i \in X$.

If $a \in L(X)$ ($a \in L^p(X)$, respectively), then partial derivatives of the element $a$ are components of the action of the universal derivation on $a$. The Jacobian conjecture for free Lie algebras over fields was proved by C. Reutenauer in [71], by V. Shpilrain in [77], and by U. U. Umirbaev in [89].

The Jacobian conjecture for free color Lie superalgebras and Lie $p$−superalgebras over fields was proved by A. A. Mikhalev in [49]. The Jacobian conjecture for free Lie algebras over rings was proved by A. A. Mikhalev and A. A. Zolotykh in [64]. The Jacobian conjecture for free commutative (non-associative) and for free anti-commutative algebras was proved by A. V. Yagzhev in [96].

The Jacobian conjecture is also true for other types of free algebras $F(X)$ (in particular, for free non-associative algebras) of finite rank as described in Section 6.2. We refer to [92] for details.

## 3.3 Primitive systems in free groups

J. Birman [7] gave a matrix characterization of automorphisms of a free group $F = F_n = \mathrm{gp}\langle x_1, \ldots, x_n \rangle$ among arbitrary endomorphisms (the "inverse function theorem") as follows. Define the matrix $J_\varphi = (d_j(y_i))_{1 \leq i,j \leq n}$ (the "Jacobian matrix" of $\varphi$), where $y_i = \varphi(x_i)$, and $d_j$ denotes partial Fox derivation (with respect to $x_j$) in the free group ring $\mathbf{Z}(F_n)$ (see Section 6.4). Then $\varphi$ is an automorphism if and only if the matrix $J_\varphi$ is invertible over $\mathbf{Z}(F_n)$.

S. Bachmuth [3] obtained an inverse function theorem for free metabelian groups on replacing the Jacobian matrix $J_\varphi$ by its image $J_\varphi^a$ over the abelianized group ring $\mathbf{Z}(F/F')$.

A. Krasnikov [38] obtained a similar result for a more general class of groups of the form $F/[R,R]$, where $R$ is an arbitrary normal subgroup of $F$.

U. Umirbaev [91] generalized Birman's result to arbitrary *primitive systems* of elements of $F_n$ (a system of elements is primitive if it can be included in a free basis):

**Theorem 15 (Umirbaev [91])** *Let $(y_1, \ldots, y_k)$, $1 \leq k \leq n$, be a system of elements of the group $F_n$. This system is primitive if and only if the matrix $(d_j(y_i))$, $1 \leq i \leq k, 1 \leq j \leq n$ is right invertible over $\mathbf{Z}(F_n)$.*

For a single primitive element, the result was previously obtained by W. Dicks and M. J. Dunwoody [18, Corollary IV.5.3].

C.K. Gupta and E. Timoshenko [30] proved a similar result for free metabelian groups, and for some other abelian-by-abelian groups. See also Romankov [72].

## 3.4 Primitive systems in free Lie (super) algebras

In what follows, $X = \{x_1, \ldots, x_n\}$; $G$ is an Abelian group; $L(X)$ is a free (color) Lie (super) algebra; $A(X)$ is a free associative algebra (it is naturally considered as an enveloping algebra for $L(X)$); $L^p(X)$ the free color Lie $p$-superalgebra. Let $L_X = L(X)$ in the case where $\operatorname{char} K = 0$, and $L_X = L^p(X)$ in the case where $\operatorname{char} K = p > 0$.

We say that a system of $G$-homogeneous elements $S \subseteq L(X)$ ($S \subseteq L^p(X)$) is *primitive* if $S$ is a subset of some set of free generators of $L(X)$ (of $L^p(X)$).

The following results are due to A. A. Mikhalev and A. A. Zolotykh [60, 61].

**Theorem 16** *Let $h$ be a $G$-homogeneous element of $L_X$. Then the following conditions are equivalent:*

(i) *the element $h$ is primitive;*

(ii) *there exist elements $m_1, \ldots, m_n \in A(X)$ such that*

$$\sum_{i=1}^n m_i \frac{\partial h}{\partial x_i} = 1.$$

Note that the statement of Theorem 16 is not valid for free Lie algebras over fields of positive characteristic. A. A. Mikhalev, U. U. Umirbaev, and A. A. Zolotykh [53] constructed the corresponding counterexample as follows. Let $K$ be a field, $\operatorname{char} K = p > 2$, $X = \{x, y, z\}$, $L(X)$ the free Lie algebra,

$$h = x + [y, z] + (\operatorname{ad} x)^p(z) \in L(X).$$

Then $h$ is not a primitive element of $L(X)$, but at the same time the rank of the left ideal of $A(X)$ generated by the elements $\frac{\partial h}{\partial x}$, $\frac{\partial h}{\partial y}$, $\frac{\partial h}{\partial z}$ is equal to 1 ($h$ is a primitive element of $L^p(X)$, but not of $L(X)$). Moreover, examples of Lie algebras over a field of prime characteristic such that these algebras are not free Lie algebras, their cohomological dimension is equal to one, and universal enveloping algebras are free associative algebras of rank two, were constructed.

**Theorem 17** *Let $h_1, \ldots, h_k$ be $G$-homogeneous elements of $L_X$.*

*Then the following conditions are equivalent:*

(i) *the system $\{h_1, \ldots, h_k\}$ is primitive;*

(ii) *there exist elements $m_{ij} \in A(X)$, $i = 1, \ldots, k$, $j = 1, \ldots, n$, such that for all $r, s = 1, \ldots, k$*

$$\sum_{i=1}^n m_{ri} \frac{\partial h_s}{\partial x_i} = \delta_{rs},$$

*i.e., $(n \times k)$-matrix $\left(\frac{\partial h_j}{\partial x_i}\right)$ is left invertible.*

Based on Theorems 16 and 17, A. A. Mikhalev and A. A. Zolotykh [61, 65, 66, 68] obtained algorithms to recognize primitive systems of elements and to construct complements of primitive systems to free generatings sets.

## 3.5 Inverse images of primitive elements

The results of this section are due to A. A. Mikhalev, V. Shpilrain, and J.-T. Yu [50]. The idea is to study the dynamics of endomorphisms of algebras $F(X)$ (those can be, in particular, free associative or Lie algebras) by considering *inverse images* of specific elements as opposed to a more traditional approach of considering *images*. In particular, it turns out that for an element $u \in F(X)$ to be an inverse image of a free generator $x_i$, is a very restrictive property.

Let $K$ be a field, $X$ a finite set, $F(X)$ a free algebra of one of the types described in Section 6.2.

**Theorem 18** *Let $u$ be an element of a free non-associative algebra $F = F(X)$, $v$ a primitive element of $F$, and $\varphi$ a monomorphism of $F$ such that $\varphi(u) = v$. Then $u$ itself is a primitive element of $F$.*

**Theorem 19** *Let $u_1, \ldots, u_k$ be elements of a free non-associative algebra $F = F(X)$ of rank $n$; $\{z_1, \ldots, z_k\}$ a primitive system of $F$, $1 \leq k \leq n$, and let $\varphi$ be a monomorphism of the algebra $F$ such that $\varphi(u_i) = z_i$, $1 \leq i \leq k$. Then $\{u_1, \ldots, u_k\}$ is a primitive system of the algebra $F$. Furthermore, if $k = n$, then $\varphi$ is an automorphism of the algebra $F$.*

**Theorem 20** *Let $U = \{u_1, \ldots, u_l\}$ be a subset of a (non-associative) algebra $F = F(X)$ which is free in one of the varieties described in Section 6.2; $\{z_1, \ldots, z_l\}$ a primitive system of $F$; $H$ the subalgebra of $F$ generated by $U$, and $\varphi$ an endomorphism of the algebra $F$ such that $\varphi(u_i) = z_i$, $1 \leq i \leq l$. Then $H$ is a retract of the algebra $F$.*

We mention that if in Theorem 20, we put $l = n$, then $\varphi$ is an automorphism, and $U$ is a set of free generators of the free algebra $F$.

For a free associative algebra of rank 2, we have:

**Theorem 21** *Let $K$ be a field, char $K = 0$, and let $a$ be an element of the free associative algebra $A_2$. Then $K\langle a \rangle$ is a retract of the algebra $A_2$ if and only if there exist a primitive element $b$ and an endomorphism $\varphi$ of the algebra $A_2$ such that $\varphi(a) = b$.*

**Theorem 22** *Let $K$ be a field, $\operatorname{char} K = 0$, $X = \{x, y\}$, $\varphi$ a monomorphism of the algebra $A_2 = K\langle X \rangle$, and let $a$ be an element of the algebra $K\langle X \rangle$ such that $\varphi(a)$ is a primitive element of $K\langle X \rangle$. Then $a$ is a primitive element of the algebra $K\langle X \rangle$.*

For a polynomial algebra, a result similar to that of Theorem 21, was proved in [85, Corollary 1.2], and a result similar to that of Theorem 22 follows from a well-known Embedding Theorem of Abhyankar and Moh [1]. For a polynomial algebra $K[x_1, \ldots, x_n]$ in more than two variables, the problem of whether or not $\varphi(u) = x_1$ for an injective $\varphi$ implies $u$ is primitive, is part of a difficult open problem known as the Abhyankar-Sathaye embedding conjecture (see [73]).

## 3.6 Free Leibniz algebras

Leibniz algebras are possible non-(anti)commutative analogs of Lie algebras. In [42, 43] these analogs were studied from the point of view of homological algebra.

Let $K$ be a field, $L$ a $K$-algebra with the bracket multiplication $[\,,\,]$. Then $L$ is a Leibniz algebra if

$$[x, [y, z]] = [[x, y], z] - [[x, z], y]$$

for all $x, y, z \in L$.

It is obvious that any Lie algebra is a Leibniz algebra, and a Leibniz algebra is a Lie algebra if the condition ($[x, x] = 0$ for all $x \in L$) is satisfied.

Let $X$ be a set, $F(X)$ the free non-associative algebra on $X$ over $K$, $I$ the two-sided ideal of $F(X)$ generated by the elements

$$[a, [b, c]] - [[a, b], c] - [[a, c], b]$$

for all $a, b, c \in F(X)$. Then the algebra $L(X) = F(X)/I$ is the free Leibniz algebra.

The free Leibniz algebra has also another presentation. Let $V$ be a vector space with the basis $X$. Then the tensor module

$$T(V) = V \oplus V^{\otimes 2} \oplus \cdots \oplus V^{\otimes n} \oplus \cdots$$

equipped with the bracket multiplication defined inductively by

$$[x, v] = x \otimes v \quad \text{for} \quad x \in T(V),\ v \in V,$$
$$[x, y \otimes v] = [x, y] \otimes v - [x \otimes v, y] \quad \text{for} \quad x, y \in T(V),\ v \in V,$$

is the free Leibniz algebra on $X$.

Let $L^l$ and $L^r$ be two copies of the Leibniz algebra $L$. We denote by $L_x$ and $R_x$ the elements of $L^l$ and $L^r$ corresponding to the universal operators of left and right multiplication on $x$. Let $I$ be the two-sided ideal of the associative tensor $K$-algebra $T(L^l \oplus L^r)$ with the identity element corresponding to the relations

$$\begin{aligned} R_{[x,y]} &= R_x R_y - R_y R_x, \\ L_{[x,y]} &= L_x R_y - R_y L_x, \\ (R_x + L_x) L_y &= 0 \text{ for any } x, y \in L. \end{aligned}$$

Then the factor algebra $UL(L) = T(L^l \oplus L^r)/I$ is the universal enveloping algebra of the Leibniz algebra $L$. The algebra $UL(L)$ is the universal multiplicative envelope of the algebra $L$ in the variety of Leibniz algebras.

Consider a mapping $d \colon L \to UL(L)$ given by $d(x) = L_x$ with $x \in L$. We set

$$I_L = \{L_x \mid x \in L\} \cdot UL(L).$$

The mapping $d \colon L \to I_L$ is the universal derivation of the algebra $L$ in the variety of Leibniz algebras. If $H$ is a subalgebra of a Leibniz algebra $L$, then we set

$$J_H = \{L_x \mid x \in H\} \cdot UL(L).$$

The following results were obtained by A. A. Mikhalev and U. U. Umirbaev [52].

**Theorem 23** *Let $x$ be an element of a Leibniz algebra $L$, $H$ a subalgebra of $L$. Then $x \in H$ if and only if $d(x) \in J_H$.*

Theorem 23 shows that the variety of all Leibniz algebras has the property of differential separability for subalgebras. This gives a negative solution of Problem 2 from [90]: the variety of all Leibniz algebras is not a Schreier variety, but has the property of differential separability for subalgebras.

The following theorem gives a solution of the Jacobian problem for free Leibniz algebras.

**Theorem 24** *Let $L$ be the free Leibniz algebra of finite rank, $\psi$ an endomorphism of $L$. Then $\psi$ is an automorphism of $L$ if and only if the Jacobian matrix $J(\psi)$ is invertible over $UL(L)$.*

An algebra $R$ is called finitely separable if for any element $a$ of $R$ and for any finitely generated subalgebra $B$ of $R$ such that $a \notin B$, there exist a finite-dimensional algebra $H$ and a homomorphism $\varphi \colon R \to H$ such that $\varphi(a) \notin \varphi(B)$.

**Theorem 25** *The free Leibniz algebras are finitely separable.*

**Corollary** *The occurrence problem for free Leibniz algebras is solvable.*

# 4 Rank of an element

In [78] and [79], V. Shpilrain introduced, for the purpose of recognizing automorphisms, a concept of *rank* of an element of a free group or a free Lie algebra as the minimal number of free generators on which an automorphic image of this element can depend. This definition naturally generalizes to arbitrary systems of elements.

Here we review some important results on the rank of a system of elements in groups and algebras.

## 4.1 Rank of free group elements

The rank of an element $u$ of a free group $F = F_n = \text{gp}\langle x_1, \ldots, x_n \rangle$ is the minimal number of generators $x_i$ on which an automorphic image of $u$ can depend. The rank of a system of elements is defined similarly.

It turns out that this rank has another interpretation.

**Theorem 26 (Umirbaev [93])** *The rank of an element $u \in F_n$ is equal to the rank of the right ideal of $\mathbf{Z}(F_n)$ generated by the elements $\left\{ \dfrac{\partial u}{\partial x_i} \,\middle|\, 1 \leq i \leq n \right\}$.*

For a system of elements, we have:

**Theorem 27 (Umirbaev [93])** *Let $Y = (y_1, \ldots, y_k)$, $1 \leq k \leq n$, be a system of elements of the group $F_n$, and let $J_Y = (d_j(y_i))$, $1 \leq i \leq k, 1 \leq j \leq n$, be the corresponding Jacobian matrix. Then the rank of $Y$ is equal to the (right) rank of the system of columns of the matrix $J_Y$, i.e., to the maximal number of right $\mathbf{Z}(F_n)$-independent columns of this matrix.*

We also note a remarkable duality here: by a result of Shpilrain [80], the (left) rank of the system of *rows* of the matrix $J_Y$ is equal to the rank (in the "usual" sense, i.e., the minimal number of generators) of the subgroup of $F_n$ generated by $Y$.

## 4.2 Rank of elements of free Lie (super) algebras

In what follows, we denote by $L_X = L(X)$ a free (color) Lie (super)algebra in the case where the field $K$ has zero characteristic, and $L_X = L^p(X)$ a free (color) Lie $p$-(super)algebra in the case where char $K = p > 2$. By $A(X)$ we denote a free associative algebra.

*The rank* of $a \in L_X$ (rank $(a)$) is the minimal number of generators from $X$ on which an element $\varphi(a)$ can depend, where $\varphi$ runs through the automorphism group of $L_X$ (in other words, rank $(a)$ is the minimal rank of a free factor of $L_X$ that contains $a$).

The following results are due to A. A. Mikhalev and A. A. Zolotykh [60, 61]. These results also give algorithms to find and to realize ranks of elements and to recognize primitive elements of free Lie (super) algebras, see [61, 65, 66, 68].

**Theorem 28** *Let $h$ be a $G$-homogeneous element of $L_X$.*

*Then* rank $(h)$ *is equal to the rank of the left ideal of $A(X)$ generated by the elements* $\left\{ \dfrac{\partial h}{\partial x} \,\middle|\, x \in X \right\}$ *(as a free left $A(X)$-module).*

A system of elements $\{h_1, \ldots, h_s\}$ of $L(X)$ (or of $L^p(X)$) has *the rank $k$* (rank $(\{h_1, \ldots, h_s\}) = k$) if $k$ is the minimal number of generators of $X$ on which the system $\{h_1, \ldots, h_s\}$ can depend, where $\varphi$ is an automorphism of $L(X)$ (of $L^p(X)$, respectively).

**Theorem 29** *Let $h_1, \ldots, h_k$ be $G$-homogeneous elements of $L_X$.*

*Then the rank of the system $\{h_1, \ldots, h_k\}$ is equal to the rank of the left $A(X)$-submodule of the $A(X)$-module*

$$A(X)^k = A(X)e_1 \oplus \cdots \oplus A(X)e_k$$

*generated by the elements*

$$\left\{ \sum_{i=1}^{k} \frac{\partial h_i}{\partial x} e_i \,\middle|\, x \in X \right\}.$$

We note that for a polynomial algebra, similar statements would not be true; the situation there is more complex. So far, we can only handle two-variable polynomial algebras; this is done in the next section.

We also note a remarkable duality here: by results of A. A. Mikhalev, V. Shpilrain, and A. A. Zolotykh, [51, 67], the (right) rank of the system of *columns* of the matrix $J_H$ is equal to the rank (in the "usual" sense, i.e., the minimal number of generators) of the subalgebra of $L_X$ generated by the set $H = \{h_1, \ldots, h_k\}$.

## 5   Two-generator algebras

We treat two-variable polynomial algebras and free associative algebras of rank two in a separate section just because there is too little known about automorphisms of polynomial or free associative algebras of bigger rank. On the other hand, there is a wealth of results in the rank 2 case, especially for two-variable polynomial algebras over $\mathbf{C}$, since powerful methods of algebraic geometry can be used in that situation.

## 5.1 Polynomial algebras in two variables

Let $P_n = K[x_1, ..., x_n]$ be the polynomial algebra in $n$ variables over a field $K$ of characteristic 0. We are going to concentrate here mainly on the algebra $P_2$.

The first description of the group $Aut(P_2)$ was given by Jung [35] back in 1942, but it was limited to the case $K = \mathbf{C}$ since he was using methods of algebraic geometry. Later on, van der Kulk [39] extended Jung's result to arbitrary ground fields. In the form we give it here, the result appears as Theorem 8.5 in P. M. Cohn's book [11].

**Theorem 30** *Every automorphism of $K[x_1, x_2]$ is a product of linear automorphisms and automorphisms of the form $x_1 \to x_1 + f(x_2); x_2 \to x_2$. More precisely, if $(g_1, g_2)$ is an automorphism of $K[x_1, x_2]$ such that $deg(g_1) \geq deg(g_2)$, say, then either $(g_1, g_2)$ is a linear automorphism, or there exists a unique $\mu \in K^*$ and a positive integer $d$ such that $deg(g_1 - \mu g_2^d) < deg(g_1)$.*

The proof given in [11] is attributed to Makar-Limanov (unpublished), with simplifications by Dicks [17].

Note that the "More precisely, ..." statement serves the algorithmic purposes: upon defining the complexity of a given pair of polynomials $(g_1, g_2)$ as the sum $deg(g_1) + deg(g_2)$, we see that Theorem 30 allows one to arrange a sequence of elementary transformations (these are linear automorphisms and automorphisms of the form $x_1 \to x_1 + f(x_2); x_2 \to x_2$) so that this complexity decreases (or, at least, does not increase) at every step, until we either get a pair of polynomials that represents a linear automorphism, or conclude that $(g_1, g_2)$ was not an automorphism of $K[x_1, x_2]$. The parallel with Nielsen's method for free groups (see Section 6.1) is obvious.

We also mention here another proof of this result (in case char $K = 0$) due to Abhyankar and Moh [1]. In fact, their method is even more similar to Nielsen's method in a free group. Many of their results are based on the following fundamental theorem which we give here only in the characteristic 0 case:

**Theorem 31** *Let $u(t), v(t) \in K[t]$ be two one-variable polynomials of degree $n \geq 1$ and $m \geq 1$. Suppose $K[t] = K[u, v]$. Then either $n$ divides $m$, or $m$ divides $n$.*

Now let's see how one can adopt a more sophisticated Whitehead's method in a polynomial algebra situation. It appears that elementary basis transformations (see Theorem 30), when applied to a polynomial $p(x_1, x_2)$, are mimicked by Gröbner transformations of a basis of the ideal of $P_2$ generated by partial derivatives of this polynomial. To be more specific, we have to give some background material first.

In the course of constructing a Gröbner basis of a given ideal of $P_n$, one uses "reductions", i.e., transformations of the following type (see [2, p.39-43]): given a pair $(p, q)$ of polynomials, set $S(p,q) = \frac{L}{l.t.(p)} p - \frac{L}{l.t.(q)} q$, where $l.t.(p)$ is the *leading term* of $p$, i.e., the *leading monomial* together with its coefficient; $L = l.c.m.(l.m.(p), l.m.(q))$ (here, as usual, $l.c.m.$ means the least common multiple, and $l.m.(p)$ denotes the leading monomial of $p$). In this section, we consider what is called "deglex ordering" in [2] - where monomials are ordered first by total degree, then lexicographically with $x_1 > x_2 > ... > x_n$.

Now a crucial observation is as follows. These Gröbner reductions appear to be of two essentially different types:

(i) *regular*, or *elementary*, transformations. These are of the form $S(p,q) = \alpha \cdot p - r \cdot q$ or $S(p,q) = \alpha \cdot q - r \cdot p$ for some polynomial $r$ and scalar $\alpha \in K^*$. This happens when the leading monomial of $p$ is divisible by the leading monomial of $q$ (or vice versa). The reason why we call these transformations *elementary* is that they can be written in the form $(p, q) \to (\alpha_1 p, \alpha_2 q) \cdot M$, where $M$ is an *elementary matrix*, i.e., a matrix which (possibly) differs from the identity matrix by a single element outside the diagonal. In case where we have more than 2 polynomials $(p_1, ..., p_k)$, we can also write $(p_1, ..., p_k) \to (\alpha_1 p_1, ..., \alpha_k p_k) \cdot M$, where $M$ is a $k \times k$ elementary matrix; elementary reduction here is actually applied to a pair of polynomials (as usual) while the other ones are kept fixed. Sometimes, it is more convenient for us to get rid of the coefficients $\alpha_i$ and write $(p_1, ..., p_k) \to (p_1, ..., p_k) \cdot M$, where $M$ belongs to the group $GE_k(P_n)$ generated by all elementary *and* diagonal matrices from $GL_k(P_n)$. It is known [87] that $GE_k(P_n) = GL_k(P_n)$ if $k \geq 3$, and $GE_2(P_n) \neq GL_2(P_n)$ if $n \geq 2$ – see [10].

(ii) *singular* transformations – these are non-regular ones.

Denote by $I_{d(p)}$ the ideal of $P_2$ generated by partial derivatives of $p$. We say that a polynomial $p \in P_n$ has a *unimodular gradient* if $I_{d(p)} = P_n$ (in particular, the ideal $I_{d(p)}$ has rank 1 in this case). Note that if the ground field $K$ is algebraically closed, then this is equivalent, by Hilbert's Nullstellensatz, to the gradient being nowhere-vanishing.

Furthermore, define the *outer rank* of a polynomial $p \in P_n$ to be the minimal number of generators $x_i$ on which an automorphic image of $p$ can depend. Note that we call it the *outer rank* here, not just the rank, to avoid confusion wit other concept(s) of rank of a polynomial.

Then we have:

**Theorem 32 (Shpilrain and Yu [83])** *Let a polynomial $p \in P_2$ have a unimodular gradient. Then the outer rank of $p$ equals 1 if and only if one can get from $(d_1(p), d_2(p))$ to $(1, 0)$ by using only elementary transformations. Or, in the matrix form: if and only if $(d_1(p), d_2(p)) \cdot M = (1,0)$ for some matrix $M \in GE_2(P_2)$.*

The proof [83] of Theorem 32 is based on a generalization of Wright's Weak Jacobian Theorem [95].

**Remark.** Elementary transformations that reduce $(d_1(p), d_2(p))$ to $(1, 0)$, can be actually chosen to be Gröbner reductions, i.e., to decrease the maximum degree of monomials *at every step* – the proof [83] is based on a recent result of Park [70].

Now we show how one can apply this result to the study of coordinate polynomials.

We call a polynomial $p \in P_n$ *coordinate* if it can be included in a generating set of cardinality $n$ of the algebra $P_n$. It is clear that the outer rank of a coordinate polynomial equals 1 (the converse is not true!). It is easy to show that a coordinate polynomial has a unimodular gradient, and again – the converse is not true! On the other hand, we have:

**Proposition ([83])** *A polynomial $p \in P_n$ is coordinate if and only if it has outer rank 1 and a unimodular gradient.*

Combining this proposition with Theorem 32 yields the following

**Theorem 33 ([83])** *A polynomial $p \in P_2$ is coordinate if and only if one can get from $(d_1(p), d_2(p))$ to $(1, 0)$ by using only elementary Gröbner reductions.*

This immediately yields an algorithm for detecting coordinate polynomials in $P_2$ (see [83]), which is similar to Whitehead's algorithm for detecting primitive elements in a free group. This algorithm is very simple and fast: it has quadratic growth with respect to the degree of a polynomial. In case $p$ is revealed to be a coordinate polynomial, the algorithm also gives a polynomial which completes $p$ to a basis of $P_2$.

It is not known whether or not there is an algorithm for detecting coordinate polynomials in $P_n$ if $n \geq 3$.

Theorems 32 and 33 also suggest the following conjecture which is relevant to an important problem known as "effective Hilbert's Nullstellensatz":

**Conjecture "G".** Let a polynomial $p \in P_2$ have a unimodular gradient. Then one can get from $(d_1(p), d_2(p))$ to $(1, 0)$ by using *at most one* singular Gröbner reduction.

**Remark.** For $n \geq 3$, Theorem 32 is no longer valid since in this case, by a result of Suslin [87], the group $GL_n(P_n) = GE_n(P_n)$ acts transitively on the set of all unimodular polynomial vectors of dimension $n$, yet there are polynomials with a unimodular gradient, but of the outer rank 2, for example, $p = x_1 + x_1^2 x_2$. The "only if" part however is valid for an arbitrary $n \geq 2$ – see [83]. It is also easy to show that one always has $\operatorname{orank}(p) \geq \operatorname{rank}(I_{d(p)})$.

Finally, we mention that our method also yields an algorithm which, given a coordinate polynomial $p \in P_2$, finds a sequence of elementary automorphisms (i.e., automorphisms of the form $x_1 \to x_1 + f(x_2)$; $x_2 \to x_2$ together with linear automorphisms) that reduces $p$ to $x_1$, and a complement of $p$ to a pair of generators of the algebra $P_2$.

## 5.2 Free associative algebras of rank two

Let $P_2 = K[x_1, x_2]$ be the polynomial algebra of rank 2 over a field $K$, and $A_2 = K\langle x_1, x_2 \rangle$ the free associative algebra of rank 2 over the same ground field.

It is well-known that the automorphism groups $Aut(P_2)$ and $Aut(A_2)$ are isomorphic, an isomorphism $Aut(A_2) \to Aut(P_2)$ being just the natural abelianization. This is due to Makar-Limanov [46] (for $K = \mathbf{C}$) and Czerniakiewicz [15] (for an arbitrary ground field). See also [11, Theorem 9.3].

In the previous section, we have described an algorithm for detecting coordinate polynomials in $P_2$. Here we use the aforementioned isomorphism between $Aut(P_2)$ and $Aut(A_2)$ to "lift" this algorithm to $A_2$ in order to detect primitive elements of the algebra $A_2$ (an element $u \in A_2$ is called *primitive* if it is an automorphic image of $x_1$; or, in other words, if there is a generating set $\{u, v\}$ of $A_2$):

**Theorem 34 (Shpilrain and Yu [84])** *There is an algorithm that distinguishes primitive elements of the algebra $A_2$ over a field of characteristic $0$.*

Here we assume that we are able to perform calculations in the ground field $K$, which basically means that, given two elements of $K$, we can decide whether or not they are equal.

Note that there is a very simple "commutator test" for deciding if a given *pair* of elements generates the algebra $A_2$ – see our Section 2.3 or [16]. The problem of distinguishing primitive elements is obviously more difficult, yet our algorithm itself is fairly simple.

Furthermore, driven by the desire to reveal *non-primitivity* of an element of $A_2$ just by inspection, we present a couple of very transparent *necessary* conditions for an element of $A_2$ to be primitive.

Denote by $J_2$ the *free special Jordan algebra* of rank 2. This is a (non-associative) unital $K$-algebra generated by the elements $x_1$ and $x_2$ of $A_2$ with respect to the binary operation $x \circ y = \frac{1}{2}(xy + yx)$. To avoid a restriction *char* $K \neq 2$, one can consider a somewhat less user-friendly definition of $J_2$ upon replacing the binary operation given above by two operations: $x \to x^2$ and $(x, y) \to xyx$.

Then we have:

**Proposition** ([84]) *For an arbitrary ground field $K$:*

(i) *The algebra $J_2$ is invariant under any automorphism of $A_2$.*

(ii) *The group $Aut(J_2)$ is isomorphic to the group $Aut(A_2)$ (and, consequently, to $Aut(P_2)$).*

**Corollary** *If $u \in A_2$ is a primitive element of $A_2$, then $u \in J_2$.*

This Corollary gives a very convenient criterion for an element of $A_2$ to be primitive. Indeed, elements of $J_2$ are characterized among the elements of $A_2$ as follows (see [8] or [33]). Define an anti-automorphism $\leftharpoonup$ of $A_2$ which re-writes every monomial backwards. For example, $(x_1 x_2)^{\leftharpoonup} = x_2 x_1$; $(x_1 x_2 x_1 x_2^2)^{\leftharpoonup} = x_2^2 x_1 x_2 x_1$ etc. Then $\leftharpoonup$ is extended to the whole $A_2$ by linearity. The elements $u \in A_2$ for which $u^{\leftharpoonup} = u$, are called *palindromic*. Then we have [8]:

– an element $u \in A_2$ belongs to $J_2$ if and only if it is palindromic.

Thus the previous corollary gives a very convenient necessary (but not sufficient) condition for primitivity:

**Corollary** *Primitive elements of $A_2$ are palindromic. (Which means, incidentally, that every homogeneous component of a primitive element is palindromic.)*

This condition is quite sensitive since the algebra $J_2$ is very small compared to the enveloping algebra $A_2$.

# 6 Background

## 6.1 The Nielsen and Whitehead methods for free groups

Let $F = F_n$ be the free group of a finite rank $n \geq 2$ with a set $X = \{x_1, ..., x_n\}$ of free generators. Let $Y = \{y_1, ..., y_m\}$ and $\widetilde{Y} = \{\widetilde{y}_1, ..., \widetilde{y}_m\}$ be arbitrary finite sets of elements of the group $F$. Consider the following elementary transformations that can be applied to $Y$:

**(N1)** $y_i$ is replaced by $y_i y_j$ or by $y_j y_i$ for some $j \neq i$ ;

**(N2)** $y_i$ is replaced by $y_i^{-1}$ ;

**(N3)** $y_i$ is replaced by some $y_j$, and at the same time $y_j$ is replaced by $y_i$.

It is understood that $y_j$ doesn't change if $j \neq i$.

One might notice that some of these transformations are redundant, i.e., are compositions of other ones. There is a reason behind that which we are going to explain a little later.

We say that two sets $Y$ and $\widetilde{Y}$ are Nielsen equivalent if one of them can be obtained from another by applying a sequence of transformations (N1)–(N3). It was proved by Nielsen that two sets $Y$ and $\widetilde{Y}$ generate the same subgroup of the group $F$ if and only if they are Nielsen equivalent. This result is now one of the central points in combinatorial group theory.

Note however that this result alone does not give an *algorithm* for deciding whether or not $Y$ and $\widetilde{Y}$ generate the same subgroup of $F$. To obtain an algorithm, we need to somehow define the *complexity* of a given set of elements, and then to show that a sequence of Nielsen transformations (N1)–(N3) can be arranged so that this complexity decreases (or, at least, does not increase) *at every step* (this is where we may need "redundant" elementary transformations!).

This was also done by Nielsen; the complexity of a given set $Y = \{y_1, ..., y_m\}$ is just the sum of the lengths of the words $y_1, ..., y_m$. We refer to [45] for details.

Nielsen's method therefore yields (in particular) an algorithm for deciding whether or not a given endomorphism of a free group of finite rank is actually an automorphism.

A somewhat more difficult problem is, given a pair of elements of a free group $F$, to find out if one of them can be taken to another by an automorphism of $F$. We call this problem *the automorphic conjugacy problem*. It was addressed by Whitehead who came up with another kind of elementary transformations in a free group:

**(W1)** For some $j$, every $x_i$, $i \neq j$, is replaced by one of the elements $x_i x_j$, $x_j^{-1} x_i$, $x_j^{-1} x_i x_j$, or $x_i$;

**(W2)** $x_i$ is replaced by $x_i^{-1}$;

**(W3)** $x_i$ is replaced by some $x_j$, and at the same time $x_j$ is replaced by $x_i$.

One might notice a similarity of the Nielsen and Whitehead transformations. However, they differ in one essential detail: Nielsen transformations are applied to arbitrary sets of elements, whereas Whitehead transformations are applied to a *fixed basis* of the group $F$.

Using (informally) matrix language, we can say that Nielsen transformations correspond to elementary rows transformations of a matrix (this correspondence can actually be made quite formal – see [80]), whereas Whitehead transformations correspond to conjugations (via changing the basis). This latter type of matrix transformation is known to be more complex, and the corresponding structural results are deeper.

There is very much the same relation between the Nielsen and Whitehead transformations in a free group.

Note also that the Whitehead transformation (W1) is somewhat more complex than its analog (N1). This is – again – to be able to arrange a sequence of elementary transformations so that the complexity of a given element (in this case, just the lexicografic length of a cyclically reduced word) would decrease (or, at least, not increase) *at every step* – see [44].

This arrangement still leaves us with a difficult problem - to find out if one of two elements *of the same complexity* (= of the same length) can be taken to another by an automorphism of $F$. This is actually the most difficult part of Whitehead's algorithm.

In one special case however this problem does not arise, namely, when one of the elements is *primitive*, i.e., is an automorphic image of $x_1$. If we have managed to reduce an element of a free group (by Whitehead transformations) to an element of length 1, we immediately conclude that it is primitive; no further analysis is needed.

Thus, the problem of distinguishing primitive elements of a free group is a relatively easy case of the automorphic conjugacy problem. As we have seen in Section 5.1, this is also the situation in a polynomial algebra.

## 6.2 Main types of free (non-associative) algebras

Let $K$ be a field, $X$ a nonempty set, $F(X)$ *the free non-associative $K$-algebra without a unit element on the set $X$*. For a subset $Z$ of $F(X)$, by $I(Z)$ we denote the two-sided ideal of $F(X)$ generated by $Z$, i.e., the smallest ideal of $F(X)$ such that $Z \subseteq I(Z)$. Let

$$Z_1 = \{ab - ba \mid a, b \in F(X)\}, \ Z_2 = \{ab + ba \mid a, b \in F(X)\}.$$

The algebra $F(X)/I(Z_1)$ is *the free commutative* (nonassociative) *algebra* on the set $X$ of free generators, and the algebra $F(X)/I(Z_2)$ is *the free anti-commutative* (non-associative) *algebra*. These algebras are the free algebras in the varieties of all commutative $K$-algebras and all anti-commutative $K$-algebras, respectively.

Let $G$ be an Abelian semigroup, $K$ a field of characteristic different from two; $\varepsilon: G \times G \to K^*$ a skew symmetric bilinear form (a commutation factor, or a bicharacter), that is,

$$\varepsilon(g, h)\, \varepsilon(h, g) = 1,$$
$$\varepsilon(g_1 + g_2, h) = \varepsilon(g_1, h)\, \varepsilon(g_2, h), \ \ \varepsilon(g, h_1 + h_2) = \varepsilon(g, h_1)\, \varepsilon(g, h_2)$$

for all $g, g_1, g_2, h, h_1, h_2 \in G$;

$$G_- = \{g \in G \mid \varepsilon(g, g) = -1\}, \quad G_+ = \{g \in G \mid \varepsilon(g, g) = +1\}.$$

A $G$-graded $K$-algebra $R = \bigoplus_{g \in G} R_g$ is *a color Lie superalgebra* if

$$[x, y] = -\varepsilon(d(x), d(y))[y, x], \quad [v, [v, v]] = 0,$$
$$[x, [y, z]] = [[x, y], z] + \varepsilon(d(x), d(y))[y, [x, z]]$$

with $d(v) \in G_-$ for $G$-homogeneous elements $x, y, z, v \in R$, where $d(a) = g$ if $a \in R_g$.

If $G = \mathbf{Z}_2$ and $\varepsilon(f, g) = (-1)^{fg}$, then color Lie superalgebra is a Lie superalgebra. If $\varepsilon \equiv 1$, then we have a $G$-graded Lie algebra (if $G = \{e\}$, then we have a Lie algebra).

We denote $(\operatorname{ad} a)(b) = [a, b]$ for all $a, b \in R$. Let $\operatorname{char} K = p > 2$. A color Lie superalgebra $R$ over $K$ is *a color Lie $p$-superalgebra* if on $G$-homogeneous components $R_g$, $g \in G_+$, we have a mapping $x \to x^{[p]}$, $d(x^{[p]}) = pd(x)$, such that for all $\alpha \in K$ and all $G$-homogeneous elements $x, y, z \in R$ with $d(x) = d(y) \in G_+$, the following conditions are satisfied:

$$(\alpha x)^{[p]} = \alpha^p x^{[p]}, \quad (\operatorname{ad}(x^{[p]}))(z) = [x^{[p]}, z] = (\operatorname{ad} x)^p(z),$$
$$(x + y)^{[p]} = x^{[p]} + y^{[p]} + \sum s_j(x, y),$$

where $js_j(x,y)$ is the coefficient on $t^{j-1}$ in the polynomial $(\mathrm{ad}\,(tx+y))^{p-1}(x)$.

If $Q$ is a $G$-graded associative algebra over $K$, then $[Q]$ denotes the color Lie superalgebra with the operation $[\,,\,]$ where $[a,b] = ab - \varepsilon(d(a),d(b))\,ba$ for $G$-homogeneous elements $a,b \in Q$.

If $\mathrm{char}\,K = p > 2$, and $x^{[p]} = x^p$ for all $G$-homogeneous elements $x$ of $Q$ with $d(x) \in G_+$, then $[Q]$ with the operation $[p]$ is a color Lie $p$-superalgebra denoted by $[Q]^p$. We consider only $G$-homogeneous elements, homomorphisms preserving the $G$-graded structure, etc.

Let $X = \{x_1,\ldots,x_n\} = \bigcup_{g\in G} X_g$ be a $G$-graded set, that is $X_g \cap X_f = \emptyset$ for $g \neq f$, $d(x) = g$ for $x \in X_g$, and let $A(X)$ be the free $G$-graded associative $K$-algebra, $L(X)$ the subalgebra of $[A(X)]$ generated by $X$. Then $L(X)$ is *the free color Lie superalgebra with the set $X$ of free generators*. In the case where $\mathrm{char}\,K = p > 2$, let $L^p(X)$ be the subalgebra of $[A(X)]^p$ generated by $X$. Then $L^p(X)$ is *the free color Lie $p$-superalgebra on $X$*.

## 6.3 The Nielsen-Schreier property for algebras

Combinatorial theory of nonassociative algebras was started by A. G. Kurosh in [40] where he proved that subalgebras of free nonassociative algebras are free.

A variety of algebras is said to be Schreier if any subalgebra of a free algebra of this variety is free. A. I. Shirshov proved in [75] that the variety of Lie algebras is Schreier (this result was also obtained by E. Witt in [94], and he also proved that the variety of all Lie $p$-algebras is Schreier). A. I. Shirshov showed in [76] that the varieties of commutative and anti-commutative algebras are Schreier.

A. A. Mikhalev in [47] and A. S. Shtern in [86] showed that the variety of Lie superalgebras is Schreier. A. A. Mikhalev proved this result for the variety of color Lie $p$–superalgebras in [48]. A well known problem is to describe all Schreier varieties of algebras. U. U. Umirbaev [90, 92] obtained new examples of Schreier varieties of algebras and gave necessary and sufficient conditions for a variety of algebras to be Schreier.

For $u \in F$, by $\ell(u) = \ell_X(u)$ we denote the degree of $u$.

A subset $M$ of $F = F(X)$ is called *independent* if $M$ is a set of free generators of the subalgebra of $F$ generated by $M$. A subset $M = \{a_i\}$ of nonzero elements of $F$ is called *reduced* if for any $i$ the leading part $a_i^\circ$ of the element $a_i$ (i.e., the sum of monomials of maximum degree) does not belong to the subalgebra of $F$ generated by the set $\{a_j^\circ \mid j \neq i\}$.

Let $S = \{s_\alpha \mid \alpha \in I\}$ be a subset of $F$. A mapping $\omega\colon S \to S' \subseteq F$ is *an elementary transformation* of $S$ if either $\omega$ is a non-degenerate linear transformation of $S$, or $\omega(s_\alpha) = s_\alpha$ for all $\alpha \in I$, $\alpha \neq \beta$, and $\omega(s_\beta) = s_\beta + f(\{s_\alpha \mid \alpha \neq \beta\})$, where $f$ is an element of a free algebra of the same variety of algebras. It is clear that elementary transformations of free generating sets induce automorphisms of the algebra $F$; such automorphisms are called *elementary*.

One can transform any finite set of elements of the algebra $F$ to a reduced set by using a finite number of elementary transformations and possibly cancelling zero elements, and the crucial point is that every reduced subset of the algebra $F$ is an independent subset (this is what is called the Nielsen property). Moreover, by using Kurosh's method, one can construct a reduced set of generators for any subalgebra of the algebra $F$. Hence any subalgebra of $F$ is free in the same variety of algebras (this is what is called the Schreier property). In [41], J. Lewin proved that for a homogeneous variety of algebras, the Nielsen and Schreier properties are equivalent. By using this equivalence, it is easy to see that automorphism groups of the corresponding free algebras (of finite rank) are generated by elementary automorphisms (for free Lie algebras, this was observed by P. M. Cohn in [9], and for free non-associative algebras by J. Lewin in [41]). We collect these facts in one statement.

**Theorem 35 ([9, 40, 41, 47, 48, 75, 76, 86, 94])** *Let $K$ be a field, $F = F(X)$ the free algebra with a finite set $X$ of free generators in one of the varieties of $K$-algebras described above. Then:*

- *Any finite subset of $F$ can be transformed into a reduced set by a finite number of elementary transformations (with possible cancellation of zeros);*

- *Any reduced subset of the algebra $F$ is an independent set;*

- *The leading part of a polynomial in a reduced set is a polynomial in leading parts of elements of this set;*

- *Every subalgebra of $F$ is free in the same variety;*

- *The automorphism group of $F$ is generated by elementary automorphisms.*

## 6.4 Free differential calculus for free groups and algebras

R. H. Fox [28] introduced free differential calculus in free group rings. Let $F = F(X)$ be the free group on a set $X$, $\mathbf{Z}F$ the integral group ring of $F$, $\Delta_F$ the augmentation ideal of $\mathbf{Z}F$, i.e., the kernel of the natural homomorphism

$$\sigma \colon \mathbf{Z}F \to \mathbf{Z}, \ \sigma\left(\sum_i (n_i f_i)\right) = \sum_i n_i,$$

$n_i \in \mathbf{Z}$, $f_i \in F$. The Fox partial derivation with respect to $x_i$ is a mapping $d_i \colon \mathbf{Z}F \to \mathbf{Z}F$ which satisfies the following conditions:

$$\begin{aligned} d_i(x_j) &= \delta_{ij}; \\ d_i(uv) &= u d_i(v) + \sigma(v)\, d_i(u); \\ d_i(ku + lv) &= k d_i(u) + l d_i(v), \end{aligned}$$

$u, v \in \mathbf{Z}F$, $k, l \in \mathbf{Z}$.

There is another interpretation of $d_i$ as follows. The ideal $\Delta_F$ is a free left $\mathbf{Z}F$-module with a free basis $\{(x_i - 1) \mid 1 \leq i \leq n\}$, and the mappings $d_i$ are projections to the corresponding free cyclic direct summands. Every element $u \in \Delta_F$ can be uniquely written in the form

$$u = \sum_{i=1}^{n} d_i(u)\,(x_i - 1).$$

For a free algebra in a variety of algebras, one can define the universal derivation of this algebra (partial derivatives of an element, the components of the universal derivation, belong to the multiplicative envelope algebra of this free algebra, see [33]).

O. G. Kharlampovich [37] seems to be the first to use (somewhat disguised though) the free differential calculus for a study of Lie algebras.

For a free associative algebra $A(X)$, denote by $A(X)^e$ the tensor product $A(X) \otimes_K A(X)$ with the multiplication given by

$$(a \otimes b)(c \otimes d) = ac \otimes db.$$

Let

$$I_A = \bigoplus_{i=1}^{|X|} (A(X)^e)_i$$

be the direct sum of $|X|$ copies of $A(X)^e$,

$$(A(X)^e)_i = A(X) \cdot dx_i \cdot A(X),\ (a \otimes b)_i = a \cdot dx_i \cdot b.$$

Let also

$$I_A = \bigoplus_{i=1}^{|X|} (A(X) \cdot dx_i \cdot A(X))_i.$$

The universal derivation of $A(X)$ is the $K$-linear mapping

$$D \colon A(X) \to I_A$$

given by $D(x_i) = 1 \cdot dx_i \cdot 1$, $D(ab) = D(a) \cdot b + a \cdot D(b)$, $D(1) = 0$. For any $a \in A(X)$, the element $D(a)$ has a unique presentation in the form $D(a) = a_1 + \cdots + a_n$, $a_i \in (A(X)^e)_i$. The component $a_i = \dfrac{\partial a}{\partial x_i}$ is called the partial derivative of $a$ with respect to $x_i$.

For example, if $X = \{x, y\}$, then

$$\frac{\partial x^2}{\partial x} = 1 \otimes x + x \otimes 1, \ \frac{\partial x^2}{\partial y} = 0;$$
$$\frac{\partial x^3}{\partial x} = 1 \otimes x^2 + x \otimes x + x^2 \otimes 1, \ \frac{\partial x^3}{\partial y} = 0.$$

If we consider the natural homomorphism of associative algebras with the identity elements $\varphi \colon A(X)^e \to K[X]$ (where $K[X]$ is the polynomial algebra), $\varphi(a \otimes b) = ab$, then the images of partial derivatives of elements of $A(X)$ are "usual" Leibniz partial derivatives.

# 7   Open problems

We start with two problems about endomorphisms of free groups.

**Problem 1** *If an endomorphism $\varphi$ of a free group of finite rank takes every primitive element to another primitive, is $\varphi$ an automorphism?*

**Problem 2** *Denote by $\operatorname{Epi}(F_n, F_k)$ the set of all homomorphisms from a free group $F_n$ onto a free group $F_k$; $n, k \geq 2$. Are there 2 elements $g_1, g_2 \in F_n$ with the following property: whenever $\varphi(g_i) = \psi(g_i), i = 1, 2$, for some homomorphisms $\varphi$, $\psi \in \operatorname{Epi}(F_n, F_k)$, it follows that $\varphi = \psi$? (In other words, every homomorphism from $\operatorname{Epi}(F_n, F_k)$ is completely determined by its values on just 2 elements).*

These two problems are also of interest when asked about various free algebras:

**Problem 3** *Let $\varphi$ be an endomorphism of a free algebra $F(X)$ of finite rank, that takes every primitive element to another primitive. Is it true that $\varphi$ is an automorphism of the algebra $F(X)$?*

We note that for free Lie algebras and color Lie superalgebras this problem was solved by A. A. Mikhalev and A. A. Zolotykh [58, 59].

**Problem 4** *Is there a polynomial $p \in P_n$ with the following property: whenever $\varphi(p) = \psi(p)$ for some non-constant-valued endomorphisms $\varphi, \psi$ of $P_n$, it follows that $\varphi = \psi$? (In other words, every non-constant-valued endomorphism of $P_n$ is completely determined by its value on just a single polynomial).*

**Problem 5** *Is it true that the intersection of two retracts of a free algebra $F(X)$ is a retract of this algebra?*

We note that Problem 5 was solved in the affirmative for free groups by G. Bergman in [6].

**Problem 6** *Construct algorithms to recognize test elements of a free algebra $F = F(X)$.*

Problem 6 was solved for free groups by L. P. Comerford [12].

**Problem 7** *Are the assertions of Theorems 28 and 29 true for a free color Lie superalgebra over the ring $\mathbf{Z}$ of integers? Are the assertions of Theorems 16 and 17 true for a free color Lie superalgebra over $\mathbf{Z}$?*

**Problem 8 (P. M. Cohn)** *Is it true that the automorphism group of a free associative algebra of finite rank is generated by elementary automorphisms? The same question for a polynomial algebra.*

Note that this problem has a positive solution for both polynomial algebra in two variables (Jung, van der Kulk [35, 39]) and the free associative algebra of rank 2 – this is due to Czerniakiewicz and Makar-Limanov [15, 46].

**Problem 9 (M. Nagata)** *Let $\varphi$ be an automorphism of the polynomial algebra $K[x, y, z]$ given by*

$$\begin{aligned} \varphi(x) &= x + z(x^2 - yz), \\ \varphi(y) &= y + 2x(x^2 - yz) + z(x^2 - yz)^2, \\ \varphi(z) &= z. \end{aligned}$$

*Is it true that this automorphism is not a composition of elementary automorphisms of $K[x, y, z]$?*

Some work on this problem has been recently done in [20].

**Problem 10 (see Section 4)** *Find criteria for an element of a polynomial or a free associative algebra to have a given rank and an algorithm to determine the rank of an element.*

**Problem 11** *Suppose $\varphi(p) = x_1$ for some monomorphism (i.e., injective endomorphism) $\varphi$ of a polynomial algebra $K[x_1, \ldots, x_n]$, $n > 2$. Is it true that $p$ is a coordinate polynomial? The same question for an element of a free associative algebra.*

When $n=2$, the affirmative answer to this problem for a polynomial algebra follows from a well-known Embedding Theorem of Abhyankar and Moh [1], and for a free associative algebra it was established in [50].

**Problem 12 (see Section 2.3)** *Characterize test polynomials of a polynomial algebra $K[x_1,\ldots,x_n]$ and of a free associative algebra $K\langle x_1,\ldots,x_n\rangle$, $n\geq 2$.*

**Problem 13** *Can any retract of $K[x_1,...,x_n]$, $n \geq 3$, be generated by algebraically independent polynomials?*

**Problem 14 (see Section 6.4)** *Is it true that an element $a$ of a free associative algebra $K\langle x,y\rangle$ is primitive if and only if there are elements $b_1, b_2$ of the algebra*

$$K\langle x,y\rangle^e = K\langle x,y\rangle \otimes_K K\langle x,y\rangle$$

*such that*

$$b_1\frac{\partial a}{\partial x} + b_2\frac{\partial a}{\partial y} = 1 \otimes 1\ ?$$

**Problem 15** *Let $a$ be an element of a free associative algebra $K\langle x,y\rangle$; $I = \mathrm{id}(a)$ the ideal of $K\langle x,y\rangle$ generated by the element $a$. Is it true that $a$ is a primitive element of $K\langle x,y\rangle$ if and only if the factor algebra $K\langle x,y\rangle/I$ is a polynomial algebra in one variable?*

# Acknowledgement

The first author thanks the Robert Black College and the Department of Mathematics, the University of Hong Kong, for their support and hospitality during his visit as a Rayson Huang Fellow when this survey was written.

# 8  Bibliography

[1] S. S. Abhyankar and T.-T. Moh, *Embeddings of the line in the plane.* J. Reine Angew. Math. **276** (1975), 148–166.

[2] W. Adams and P. Loustaunau, *An introduction to Gröbner bases.* American Mathematical Society, Providence, 1994.

[3] S. Bachmuth, *Automorphisms of free metabelian groups*, Trans. Amer. Math. Soc. **118** (1965), 93-104.

[4] Yu. A. Bahturin, A. A. Mikhalev, M. V. Zaicev, and V. M. Petrogradsky, *Infinite Dimensional Lie Superalgebras*. Walter de Gruyter Publ., Berlin–New York, 1992.

[5] H. Bass, E. Connell, and D. Wright, *The Jacobian conjecture: reduction on degree and formal expansion of the inverse*. Bull. Amer. Math. Soc. **7** (1982), 287–330.

[6] G. Bergman, *Supports of derivations, free factorizations, and ranks of fixed subgroups in free groups*. Trans. Amer. Math. Soc., in press.

[7] J. S. Birman, *An inverse function theorem for free groups*. Proc. Amer. Math. Soc. **41** (1973), 634–638

[8] P. M. Cohn, *On homomorphic images of special Jordan algebras*, Canadian J. Math. **6** (1954), 253–264.

[9] P. M. Cohn, *Subalgebras of free associative algebras*. Proc. London Math. Soc. (3) **14** (1964), 618–632.

[10] P. M. Cohn, *On the structure of the $GL_2$ of a ring*. Inst. Hautes Études Sci. Publ. Math. **30** (1966), 365–413.

[11] P. M. Cohn, *Free Rings and Their Relations. 2nd Ed.*. Academic Press, London, 1985.

[12] L. P. Comerford, *Generic elements of free groups*. Arch. Math. (Basel) **65** (1995), no. 3, 185–195.

[13] E. Connell and J. Zweibel, *Subrings invariant under polynomial maps*. Houston J. Math. **20** (1994), 175–185.

[14] D. Costa, *Retracts of polynomial rings*. J. Algebra **44** (1977), 492–502.

[15] A. J. Czerniakiewicz, *Automorphisms of a free associative algebra of rank 2*, I, II. Trans. Amer. Math. Soc. **160** (1971), 393–401; **171** (1972), 309–315.

[16] W. Dicks, *A commutator test for two elements to generate the free algebra of rank two*. Bull. London Math. Soc. **14** (1982), 48–51.

[17] W. Dicks, *Automorphisms of the polynomial ring in two variables*. Publ. Sec. Mat. Univ. Autonoma Barcelona **27** (1983), 155–162.

[18] W. Dicks and M. J. Dunwoody, *Groups acting on graphs*. Cambridge University Press, 1989.

[19] W. Dicks and J. Lewin, *A Jacobian conjecture for free associative algebras*. Comm. Algebra **10** (1982), 1285–1306.

[20] V. Drensky, J. Gutierrez, and J.-T. Yu, *Gröbner bases and the Nagata automorphism*. J. Pure Appl. Algebra, to appear.

[21] V. Drensky and J.-T. Yu, *Test polynomials for automorphisms of polynomial and free associative algebras*. Preprint, 1996.

[22] V. Drensky and J.-T. Yu, *Orbits in polynomial and relatively free algebras of rank two*. Commun. Algebra, to appear.

[23] V. G. Durnev, *The Maltsev-Nielsen equation in a free metabelian group of rank two*. Math. Notes **46** (1989), 927–929.

[24] A. van den Essen, *Seven lectures on polynomial automorphisms*. In A. van den Essen (Ed.), *Automorphisms of Affine Spaces*, Kluwer Academic Publishers, Dordrecht, 1995, 3–39.

[25] A. van den Essen and V. Shpilrain, *Some combinatorial questions about polynomial mappings*. J. Pure Appl. Algebra **119** (1997), 47–52.

[26] B. Fine, N. Isermann, G. Rosenberger, and D. Spellman, *Test words, generic elements and almost primitivity*. Preprint, 1996.

[27] E. Formanek, *Observations about the Jacobian conjecture*. Houston J. Math. **20** (1994), 369–380.

[28] R. H. Fox, *Free differential calculus. I. Derivations in free group rings*. Ann. of Math. (2) **57** (1953), 547–560.

[29] N. Gupta and V. Shpilrain, *Nielsen's commutator test for two-generator groups*. Math. Proc. Cambridge Phil. Soc. **114** (1993), 295–301.

[30] C. K. Gupta and E. I. Timoshenko, *Primitivity in the free groups of the variety $\mathcal{A}_m\mathcal{A}_n$*. Commun. Algebra **24** (1996), 2859–2876.

[31] S. V. Ivanov, *On certain elements of free groups*. J. Algebra, to appear.

[32] S. V. Ivanov, *On endomorphisms of free groups that preserve primitivity*. Preprint.

[33] N. Jacobson, *Structure and Representations of Jordan Algebras*. Amer. Math. Soc., Providence, Rhode Island, 1968.

[34] Z. Jelonek, *A solution of the problem of van den Essen and Shpilrain*. Preprint IMUJ 1996/19, Krakow.

[35] H. W. E. Jung, *Über ganze birationale Transformationen der Ebene.* J. reine angew Math. **184** (1942), 161–174.

[36] O. Keller, *Ganze Cremona-Transformationen.* Monatsh. Math. Phys. **47** (1939), 299–306.

[37] O. G. Kharlampovich, *Lyndon condition for solvable Lie algebras.* Izv. Vyssh. Uchebn. Zaved. Mat. **1984**, no. 9, 50–59.

[38] A. F. Krasnikov, *Generators of the group $F/[N,N]$.* Math. Notes **24** (1979), 591–594.

[39] W. van der Kulk, *On polynomial rings in two variables.* Nieuw Archief voor Wisk. (3) **1** (1953), 33–41.

[40] A. G. Kurosh, *Nonassociative free algebras and free products of algebras.* Mat. Sb. **20** (1947), 239–262.

[41] J. Lewin, *On Schreier varieties of linear algebras.* Trans. Amer. Math. Soc. **132** (1968), 553–562

[42] J.-L. Loday, *Une version non commutative des algèbres de Lie: les algèbres de Leibniz.* L'Enseignement Math. **39** (1993), 269–293.

[43] J.-L. Loday and T. Pirashvili, *Universal enveloping algebras of Leibniz algebras and (co)homology.* Math. Ann. **296** (1993), 139–158.

[44] R. Lyndon and P. Schupp, *Combinatorial Group Theory.* Series of Modern Studies in Math. **89**. Springer-Verlag, 1977.

[45] W. Magnus, A. Karrass, and D. Solitar, *Combinatorial Group Theory.* Wiley, New York, 1966.

[46] L. G. Makar-Limanov, *The automorphisms of the free algebra of two generators.* Funktsional. Anal. i Prilozhen. **4** (1970), no. 3, 107–108.

[47] A. A. Mikhalev, *Subalgebras of free color Lie superalgebras.* Mat. Zametki **37** (1985), No. 5, 653–661. English translation: Math. Notes, **37** (1985), 356–360.

[48] A. A. Mikhalev, *Subalgebras of free Lie p-superalgebras.* Mat. Zametki **43** (1988), No. 2, 178–191. English translation: Math. Notes **43** (1988), 99–106.

[49] A. A. Mikhalev, *On right ideals of the free associative algebra generated by free color Lie (p-)superalgebras.* Uspekhi Mat. Nauk **47** (1992), no. 5, 187–188. English translation: Russian Math. Surveys. **47** (1992), No. 5, 196–197.

[50] A. A. Mikhalev, V. Shpilrain, and J.-T. Yu, *On the dinamics of homomorphisms of free algebras*. Preprint, 1997.

[51] A. A. Mikhalev, V. Shpilrain, and A. A. Zolotykh, *Subalgebras of free algebras*. Proc. Amer. Math. Soc. **124** (1996), 1977–1983.

[52] A. A. Mikhalev and U. U. Umirbaev, *Subalgebras of free Leibniz algebras*. Commun. Algebra, in press.

[53] A. A. Mikhalev, U. U. Umirbaev, and A. A. Zolotykh, *An example of a nonfree Lie algebra of cohomological dimension one*. Uspekhi Mat. Nauk **49** (1994), no. 1, 203–204. English translation: Russian Math. Surveys **49** (1994), no. 1, p. 254.

[54] A. A. Mikhalev and J.-T. Yu, *Test elements and retracts of free Lie algebras*. Commun. Algebra **25** (1997), 3283–3289.

[55] A. A. Mikhalev and J.-T. Yu, *Test elements, retracts and automorphic orbits of free algebras*, Internat. J. Algebra Comput., in press.

[56] A. A. Mikhalev and J.-T. Yu, *Automorphic orbits in free algebras of rank two*. Preprint, 1997.

[57] A. A. Mikhalev, J.-T. Yu, and A. A. Zolotykh, *Images of coordinate polynomials*. Algebra Colloquium **4** (1997), 159–162.

[58] A. A. Mikhalev and A. A. Zolotykh, *Endomorphisms of free Lie algebras preserving primitivity of elements are automorphisms*. Uspekhi Matem. Nauk **48** (1993), No. 6, 149–150. English translation: Russian Math. Surveys **48** (1993), No. 6, 189–190.

[59] A. A. Mikhalev and A. A. Zolotykh, *Automorphisms and primitive elements of free Lie superalgebras*. Commun. Algebra **22** (1994), 5889–5901.

[60] A. A. Mikhalev and A. A. Zolotykh, *The rank of an element of the free color Lie (p-)superalgebra*. Doklady Akad. Nauk **334** (1994), 690–693. English translation: Russian Acad. Sci. Dokl. Math. **49** (1994), No. 1, 189–193.

[61] A. A. Mikhalev and A. A. Zolotykh, *Rank and primitivity of elements of free color Lie (p-)superalgebras*. Intern. J. Algebra Comput. **4** (1994), 617–656.

[62] A. A. Mikhalev and A. A. Zolotykh, *Test elements for monomorphisms of free Lie algebras and superalgebras*. Commun. Algebra **23** (1995), 4995–5001.

[63] A. A. Mikhalev and A. A. Zolotykh, *Endomorphisms of free associative algebras over commutative rings and their Jacobian matrices*. Fundamental'naya i Prikladnaya Matematika **1** (1995), 177–190.

[64] A. A. Mikhalev and A. A. Zolotykh, *An inverse function theorem for free Lie algebras over commutative rings*. Algebra Colloquium **2** (1995), no. 3, 213–220.

[65] A. A. Mikhalev and A. A. Zolotykh, *Combinatorial Aspects of Lie Superalgebras*. CRC Press, Boca Raton, New York, 1995.

[66] A. A. Mikhalev and A. A. Zolotykh, *Algorithms for primitive elements of free Lie algebras and superalgebras*. Proc. ISSAC-96, ACM Press, New York, 1996, 161–169.

[67] A. A. Mikhalev and A. A. Zolotykh, *Ranks of subalgebras of free Lie superalgebras*. Vestnik Mosk. Univ. Ser.1. Mat. Mekh., **1996**, no. 2, 36–40. English translation: in Moscow Univ. Math. Bull.

[68] A. A. Mikhalev and A. A. Zolotykh, *Complex of algorithms for computation in Lie superalgebras*. Programmirovanie **1997**, no. 1, 12–23.

[69] J. Nielsen *Die Isomorphismen der allgemeinen, unendlichen Gruppe mit zwei Erzeugenden*. Math. Ann. **78** (1918), 385–397.

[70] H. Park, *A Computational Theory of Laurent Polynomial Rings and Multi-dimensional FIR Systems*. Ph.D thesis, University of California at Berkeley, 1995.

[71] C. Reutenauer, *Applications of a noncommutative Jacobian matrix*. J. Pure Appl. Algebra **77** (1992), 169–181.

[72] V. A. Romankov, *Criteria for the primitivity of a system of elements of a free metabelian group*. Ukrainskii Mat. Zh. **43** (1991), 996–1002.

[73] A. Sathaye, *On linear planes*. Proc. Amer. Math. Soc. **56** (1976), 1–7.

[74] A. H. Schofield, *Representations of Rings over Skew Fields*. London Math. Soc. Lecture Note Ser. **92** (1985).

[75] A. I. Shirshov, *Subalgebras of free Lie algebras*. Mat. Sb. **33** (1953), 441–452.

[76] A. I. Shirshov, *Subalgebras of free commutative and free anti-commutative algebras*. Mat. Sb. **34** (1954), 81–88.

[77] V. Shpilrain, *On generators of $L/R^2$ Lie algebras*. Proc. Amer. Math. Soc. **119** (1993), 1039–1043.

[78] V. Shpilrain, *Recognizing automorphisms of the free groups*. Arch. Math. **62** (1994), 385–392.

[79] V. Shpilrain, *On the rank of an element of a free Lie algebra*. Proc. Amer. Math. Soc. **123** (1995), 1303–1307.

[80] V. Shpilrain, *On monomorphisms of free groups*. Arch. Math. **64** (1995), 465–470.

[81] V. Shpilrain, *Test elements for endomorphisms of the free groups and algebras*. Israel J. Math. **92** (1995), 307–316.

[82] V. Shpilrain, *Generalized primitive elements of a free group*. Preprint.

[83] V. Shpilrain and J.-T. Yu, *Polynomial automorphisms and Gröbner reductions*. J. Algebra, to appear.

[84] V. Shpilrain and J.-T. Yu, *On generators of polynomial algebras in two commuting or non-commuting variables*. J. Pure Appl. Algebra, to appear.

[85] V. Shpilrain and J.-T. Yu, *Polynomial retracts and the Jacobian conjecture*. Trans. Amer. Math. Soc., to appear.

[86] A. S. Shtern, *Free Lie superalgebras*. Sibirsk. Mat. Zh. **27** (1986), 170–174. English translation: in Sib. Math. J.

[87] A. A. Suslin, *On the structure of the special linear group over polynomial rings*. Math. USSR Izv. **11** (1977), 221–238.

[88] E. C. Turner, *Test words for automorphisms of free groups*. Bull. London Math. Soc. **28** (1996), 255–263.

[89] U. U. Umirbaev, *Partial derivatives and endomorphisms of some relatively free Lie algebras*. Sibirsk. Mat. Zh. **34** (1993), no. 6, 179–188. English translation: in Sib. Math. J.

[90] U. U. Umirbaev, *On Schreier varieties of algebras*. Algebra i Logika **33** (1994), 317–340. English translation: in Algebra and Logic.

[91] U. U. Umirbaev, *Primitive elements of free groups*. Uspekhi Mat. Nauk **49** (1994), no. 2, 175–176. English transalation: in Russian Math. Surveys.

[92] U. U. Umirbaev, *Universal derivations and subalgebras of free algebras*. In *Proc. 3rd Int. Conf. in Algebra (Krasnoyarsk, 1993)*. Walter de Gruyter, 1996, 255–271.

[93] U. U. Umirbaev, *On ranks of elements of free groups*. Fundamentalnaya i Prikladnaya Matematika **2** (1996), 313–315.

[94] E. Witt, *Die Unterringe der freien Lieschen Ringe*. Math. Z. **64** (1956), 195–216.

[95] D. Wright, *The amalgamated free product structure of $GL_2(k[X_1, ..., X_n])$ and the weak Jacobian theorem for two variables*. J. Pure Appl. Algebra **12** (1978), 235–251.

[96] A. V. Yagzhev, *Endomorphisms of free algebras*. Sib. Math. J. **21** (1980), no. 1, 133–141.

[97] J.-T. Yu, *On the Jacobian conjecture: reduction of coefficients*. J. Algebra **171** (1995), 213–219.

Alexander A. Mikhalev: Department of Mechanics and Mathematics, Moscow State University, Moscow 119899, Russia
*e-mail address:* aamikh@cnit.math.msu.su

Vladimir Shpilrain: Department of Mathematics, City College of New York, New York, NY 10031
*e-mail address:* shpil@groups.sci.ccny.cuny.edu

Jie-Tai Yu: Department of Mathematics, The University of Hong Kong, Pokfulam Road, Hong Kong
*e-mail address:* yujt@hkusua.hku.hk